\documentclass[11pt]{article}

\usepackage{inputenc}
\usepackage{makeidx}
\usepackage{amssymb}
\usepackage{amsfonts}
\usepackage{amsmath}
\usepackage{euscript}
\usepackage{theorem}
\usepackage{enumerate}
\usepackage{graphics}

\theorembodyfont{\rm}      

\theoremstyle{change}      

\begingroup

\endgroup

\begingroup

\endgroup

\begingroup

\endgroup

\begin{document}

$$\text{\bf\Large  A Closed Formula for the Product}$$
$$\text{\bf\Large in Simple Integral Extensions}$$


\begin{center}{{\large {\sc Natalio H. Guersenzvaig}}\\
Av Corrientes 3985 6A, (1194) Buenos Aires, Argentina\\
email: nguersenz@fibertel.com.ar}
\end{center}

\begin{center}
and
\end{center}

\begin{center}{{\large {\sc Fernando Szechtman\begin{footnote}{Corresponding author}\end{footnote}}}\\
Department of Mathematics and Statistics, University of Regina, Saskatchewan, Canada\\
email: fernando.szechtman@gmail.com}
\end{center}


\begin{abstract}
Let $\xi$ be an algebraic number and let $\alpha,\beta\in \mathbb
Q[\xi]$. An explicit formula for the coordinates of the product
$\alpha\beta$ is given in terms of the coordinates of $\alpha$ and
$\beta$ and the companion matrix of the minimal polynomial of
$\xi$. The formula as well as its proof extend to fairly general
simple integral extensions.
\end{abstract}


\noindent{\bf Keywords:} Integral extension, companion matrix.

\smallskip

\noindent{\bf AMS Classification:} 15A24

\bigskip

Let $\xi$ be an algebraic number of degree $n$ and minimal
polynomial $f\in \mathbb Q[X]$, say
$f=X^{n}+a_{n-1}X^{n-1}+\cdots+a_1X+a_0$. Given $\alpha\in \mathbb
Q(\xi)$ let $[\alpha]$ stand for the coordinates of $\alpha$
relative to the basis $1, \xi, \dots , \xi^{n-1}$. For $g,h$ in
$\mathbb Q[X]$ let $\beta=g(\xi),\gamma=h(\xi)\in\mathbb Q(\xi)$.
Obviously $[\beta+\gamma]=[\beta]+[\gamma]$. On the other hand
$[\beta\gamma]$ is usually obtained  either by reducing high
powers of $\xi$ to smaller ones using the relation $\xi^n=
-a_{n-1}\xi^{n-1} - \cdots - a_1\xi -a_0$, or by finding the
remainder of dividing $gh$ by $f$. We wonder if there is also a
closed formula for $[\beta\gamma]$ in terms of $[\beta]$
and~$[\gamma]$ that replaces these procedures. Surprisingly, such
a formula exists, and it is so natural that it works in greater
generality, as described below.

\medskip
\noindent{\bf Theorem 1.} {\sl Let $R$ be a ring with identity.
Assume that $\xi$ belongs to an overing of $R$, that $\xi$
commutes with all elements of $R$ and that
 $f=X^{n}+a_{n-1}X^{n-1}+\cdots+a_1X+a_0\in R[X]$ is the minimal polynomial of $\xi$ over $R$,  i.e., that $\xi$ is a root of $f$ but not a root of any polynomial in $R[X]$ of degree less than $n$. Let
$\alpha,\beta\in R[\xi]$ have coordinates $[\alpha],[\beta]$
relative to the basis $1,\xi,\dots,\xi^{n-1}$. Then
\begin{equation*}
\label{form1}
[\alpha\beta]=(I\; C\;\dots\; C^{n-1})([\alpha]\otimes[\beta]),
\end{equation*}
where $C\in M_n(R)$ is the companion matrix of $f$,
$$
C=\left(%
\begin{array}{ccccc}
  0 & 0 & \cdots & 0 & -a_0 \\
  1 & 0 & \cdots & 0 & -a_1 \\
  0 & 1 & \cdots & 0 & -a_2 \\
  \vdots & \vdots & \cdots & \vdots & \vdots \\
  0 & 0 & \cdots & 1 & -a_{n-1} \\
\end{array}%
\right),
$$
and for $x,y\in R^n$ their left Kronecker product $x\otimes y\in
R^{n^2}$ is given by
$$
x\otimes y=\left(%
\begin{array}{c}
  xy_1 \\
  \vdots \\
  xy_n \\
\end{array}%
\right).
$$}
\noindent{\sl Proof.} Observe first of all that the coefficients
of $f$ belong to the center of $R$.  Indeed, let $y\in R$. Then
$0=yf(\xi)-f(\xi)y=(ya_0-a_0y)+(ya_1-a_1y)\xi+\cdots+(ya_{n-1}-a_{n-1}y)\xi^{n-1}$.
Thus $R[C]$ is a ring with $C$ in its center. Clearly $f(C)=0$ and
$I, C, \dots , C^{n-1}$ is a basis of the $R$-module $R[C]$. Let
$[A]$ be the coordinates of $A\in R[C]$ in this basis. The map
$R[C]\to R[\xi]$ given by $p(C)\mapsto p(\xi)$, $p\in R[X]$, is a
ring isomorphism that preserves coordinates, so it suffices to
prove the result in $R[C]$. To this end, note that if $A\in R[C]$
\begin{equation}
\label{uno} A = ([A] \;C[A] \;\dots \;C^{n-1}[A]).
\end{equation}
Indeed, let $e_1,\dots,e_n$ be the canonical basis of $R^n$. We
have $A = y_0I + y_1C + \cdots + y_{n-1}C^{n-1}$ with $y_j\in R$,
so $Ae_1 = y_0e_1 + y_1e_2 + \cdots + y_{n-1}e_n = [A]$. But $Ae_j
= AC^{j-1}e_1=~C^{j-1}Ae_1$, so $Ae_j=C^{j-1}[A]$, for all $2\le
j\le n$. Thus the matrices in question have the same columns.

We know from (\ref{uno}) that $[AB]$ is the first column of $AB$,
i.e. $[AB] = A[B]$, also by  (\ref{uno}). Let 0 stand for the zero
column vector of length $n$. Applying  (\ref{uno}) once more gives
\begin{align*}
\quad\quad\quad [AB] & = ([A] \;C[A]\; \dots \; C^{n-1}[A])[B]\\
& = (I\; C \dots C^{n-1})\!\left(%
\begin{array}{cccc}
  [A] & 0 & \cdots & 0 \\
  0 & [A] & \cdots & 0 \\
  \vdots & \vdots & \cdots & \vdots \\
  0 & 0 & \cdots & [A] \\
\end{array}%
\right)\![B]\\
& = (I\; C\;\dots\; C^{n-1})
([A]\otimes[B]).\qquad\qquad \qquad\blacksquare
\end{align*}

As a typical example, let $R=\mathbb R$ and
$\xi=i$, so that $f(X)=X^2+1$. Then
$$
[(b_0+b_1 \xi)(c_0+c_1 \xi)]=\left(%
\begin{array}{ccccc}
  1 & 0 & | & 0 & -1 \\
  0 & 1 & | & 1 & 0 \\
\end{array}%
\right)\left(%
\begin{array}{c}
  b_0c_0 \\
  b_1c_0 \\
  b_0c_1 \\
  b_1c_1 \\
\end{array}%
\right)=\left(%
\begin{array}{c}
  b_0c_0-b_1c_1 \\
 b_1c_0+b_0c_1 \\
\end{array}%
\right),
$$
as expected from the usual definition  of the product in $\mathbb
R^2$ that produces $\mathbb C$. This can  be generalized by
defining  an operation on the column space $R^n$ that makes it a
ring isomorphic to $R[\xi]$ via the map $R[\xi]\mapsto R^n$ given
 by $\alpha\mapsto [\alpha]$. We just define
\begin{equation*}
\label{g8} x\cdot y=(I\; C\;\cdots \; C^{n-1})(x\otimes y),\quad
x,y\in R^n.
\end{equation*}

Suppose next $A\in M_n(R)$ has coefficients in the center of $R$.
Then its characteristic polynomial, say $f$, is monic of degree
$n$ with all its coefficients in the center of $R$ and satisfies $f(A)=0$. However, $A$ need not have a minimal
polynomial, and even if it has one it may be harder to determine
than~$f$. For cases such as these we have the result below, where
$R_n[X]=\{p\in R[X]: \,p=0 \text{ or } \deg \,p <n\}$ and $[p]$
denotes the column vector in $R^n$ formed by the coefficients of
$p\in R_n[X]$.

\medskip
\noindent{\bf Theorem 2.} {\sl Let $R$ be a ring with identity.
Assume that $\xi$ belongs to an overing of $R$, that~$\xi$
commutes with all elements of $R$, and that $f(\xi)=0$ for some
monic polynomial $f\in R[X]$ of degree $n$ with all its coefficients in the center of $R$.
Let $C$ denote the companion matrix of $f$. Then for all $g, h\in R_n[X]$ we have
$$g(\xi)h(\xi) = (g\odot h)(\xi),$$ where $g\odot h$ denotes the polynomial in
$R_n[X]$ whose coefficients are given by
\begin{equation*}
 [g\odot h]= (I\; C\;\cdots\;
C^{n-1})([g]\otimes[h]).
\end{equation*}}
\noindent{\sl Proof.} This follows from Theorem 1 via the
epimorphism $R[C]\mapsto R[\xi]$.$\quad\qquad \qquad\qquad
\blacksquare$

\medskip
As an example, let $\xi =A \in M_3(R)$ have coefficients in the
center of $R$ and characteristic polynomial $f=X^3-1$. Let
$g=b_0+b_1X+b_2X^2$, $h=c_0+c_1X+c_2X^2$ be in $R_3[X]$. Then
$g(A)h(A)= (g\odot h)(A)$, where
\begin{align*}
\label{example}
[g\odot h] &=\left(\begin{array}{ccccccccccc}
  1 & 0 & 0 &|& 0 & 0 & 1& |&0&1&0\\
  0 & 1 & 0 &|& 1 & 0 & 0 &|&0&0&1\\
  0 & 0 & 1&| & 0 & 1 & 0&|&1&0&0\end{array}\right)
 \!\!\left(\begin{array}{c}
  b_0c_0 \\
  \vdots \\
  b_1c_2\\
  b_2c_2\\
\end{array}\right)= \left(\begin{array}{c}
b_0c_0+b_2c_1 +b_1c_2\\
b_1c_0 +b_0c_1+b_2c_2 \\
b_2c_0+ b_1c_1+ b_0c_2\end{array}\right).
\end{align*}
This identity can be easily generalized for $f=X^n-1$ using the formula $[g\odot h] = g(C)[h]$, which is particularly useful when $g(C)$ has a specified form, as in this example, related to the well known {\sl circulant} matrices.

\end{document}